# The Collapse of the Hilbert Program:
# A Variation on the Gödelian Theme[*]

Saul A. Kripke

*Abstract*. The Hilbert program was actually a specific approach for proving consistency. Quantifiers were supposed to be replaced by ε-terms. ε*xA(x)* was supposed to denote a witness to $\exists x A(x)$, arbitrary if there is none. The Hilbertians claimed that in any proof in a number-theoretic system *S*, each ε-term can be replaced by a numeral, making each line provable and true. This implies that *S* must not only be consistent, but also 1-consistent ($\Sigma_1^0$-correct).

Here we show that if the result is supposed to be provable within *S*, a statement about all $\Pi_2^0$ statements that subsumes itself within its own scope must be provable, yielding a contradiction. The result resembles Gödel's but arises naturally out of the Hilbert program itself.

In contrast to what seems to be the case today, in the early part of the 20[th] century many leading mathematicians were very concerned with problems in the foundations of mathematics. They thought there was a great problem awaiting their contributions. Hilbert was no exception. The main idea of what has become known as the Hilbert program was of course to prove the consistency of various systems of mathematics by 'finitary' means. He and his followers thought that this would show that a very narrow class of statements – those that they regarded as having a finitary meaning

---

[*] This paper is based on a transcript of a lecture given at Indiana University on October 15, 2007 in honor of the inauguration of President Michael A. McRobbie. An abstract of another version of this talk (given at the *2008 Winter Meeting of the Association for Symbolic Logic* on December 30[th], 2008) was published in *The Bulletin of Symbolic Logic* 15(2), 2009: 229-231.

I am indebted to Burton Dreben for his insistence that the Hilbert program or approach (*Ansatz*) was not merely to prove the consistency of mathematics by finitary means, but was a specific program for interpreting proofs. Thus, as Dreben emphasized, it is a kind of constructive model theory.

– would therefore be shown to be provable by finitary 'metamathematical' means, even though they were proved in systems that we must consider too 'infinitistic' to regard as true.

The antagonism between Hilbert and Brouwer in much of the foundational debate of this period is famous.[1] Yet, it is worth remarking that Hilbert's actual program in the 1920s advocated an official 'finitism' which, though never fully defined, seems to have found valid a much *narrower* class of arguments than those admitted by Brouwer. The class of statements that had a clear finitary meaning was also narrow. (Hilbert had little interest in the mathematical statements and arguments accepted by Brouwer.) A system containing infinitary mathematics was supposed to be justified by a finitary metamathematical argument that would show that all statements provable in the system in the narrowly meaningful class can be proved without resort to the dubious statements used in the proof.

Now, Hilbert's interests are only partially described by the formulation of the Hilbert program in the 1920s. They actually progressed to that position. In the preceding years, for example, he was interested in the approach of *Principia Mathematica* and wished to have Russell come to Göttingen (see Sieg 1999 for a history).[2] But although he oscillated in his views, he really seemed to have hankered after a position in which classical mathematics will really be true, even though officially it largely consists only of 'ideal' statements used to prove 'contentual' statements. And his famous slogan, "No one shall be able to drive us from the Paradise that Cantor created for us" in his relatively late (and rather obscure in his claim to have solved the continuum problem) paper "On the Infinite" (1926: 376), seems to confirm that this was a long-held position.

We all know that the story with Hilbert's program doesn't have a happy ending, but we should also remember that it did not 'collapse' (as I say in my title) without a legacy – both in the notion of a formal system as rigorously formulated and independent of its interpretation, and in the subject of proof theory (in which a lot of Hilbert's standards have had to be relaxed to carry some goals of the program through). Unlike *Principia Mathematica*, Hilbert and his followers rigorously separated the purely formal, syntactic, formulation of a system, from its interpretation. As Gödel has remarked, *Principia* was a considerable backwards step from Frege's work in the rigor of its formalization.[3] Frege's standards of formalization were in fact the modern ones, even though (as far as I know) he did not emphasize what the standards should be. Hilbert and his school did so.

---

[1] For the culmination of the *Grundlagenstreit* between Hilbert and Brouwer, see van Dalen (1990). His paper describes Hilbert's ultimately successful, though legally unsupportable, effort to have Brouwer removed from the editorial board of the journal *Mathematische Annalen*, the leading mathematics journal of his time. He describes the reaction of various members of the editorial board, usually major mathematicians, although Albert Einstein wasn't a mathematician and was a member of the board. The title of van Dalen's paper ("The War of the Frogs and the Mice") is taken from an Ancient Greek parody of the *Iliad*, and it is often used to describe a vehement, but essentially trivial dispute. Einstein once used the phrase to refer to this dispute. It left Brouwer very depressed and temporarily unproductive. van Dalen also relates that Hilbert and Brouwer were originally on very good terms.

[2] Sieg (1999) states that Hilbert's attitude really narrowed and became more defensive as time went on – originally, perhaps, he had hoped simply to justify the absolute meaningfulness and correctness of ordinary mathematics (say, in *Principia Mathematica*) outright without revision, whereas eventually he settled for mounting only a very conservative defense.

Hilbert and Ackermann's famous textbook (1928) formulates, among other things, first-order logic without identity (really already known to Hilbert much earlier); and the problems of completeness and decidability with respect to validity in any non-empty domain are clearly stated. In other lectures Hilbert did include identity as a part of first-order logic. Gödel (1930) proved completeness for both systems without and with identity.[4] Regarding proof theory, Hilbert formulated an approach to the subject, as well as the subject itself. This subject is still alive today, and is even related to some current problems in computer science and other formal disciplines. But *complete* realization of Hilbert's original program of the 1920s is usually thought to have been shown to be impossible.[5]

The details of Hilbert's approach, and the reasons it convinced a generation of logicians to believe that it obviously must succeed and that only technical work remained to finish it, are generally left to more specialized textbooks. As we will see, the 'Hilbertians' had ample reason to assume that it would succeed. Indeed, Gödel's incompleteness theorems came as a shock to them – especially the second incompleteness theorem, in which he gave a rigorous proof that a reasonably strong consistent formal system cannot prove its own consistency. Gödel himself wrote that his work did not disprove the viability of the Hilbert program. Perhaps there are finitary

---

[3] Gödel says: "It is to be regretted that it (*Principia*) is so lacking in formal precision in the foundations that it represents in this respect a considerable step backwards as compared with Frege" (1944: 120).

My own personal view is that *Principia* should not simply be viewed as an unsuccessful attempt to give a Hilbertian formal system, but it is some kind of axiomatic theory of propositions and 'propositional functions'. This is all I can say here.

[4] I am indebted here to correspondence with Richard Zach. See also Ewald (2019), especially §8.

[5] Let me emphasize that Paul Bernays and Wilhelm Ackermann were crucial to the development of Hilbertian logic and proof theory as well as the program itself. Bernays actually wrote the book *Grundlagen der Mathematik*, by Hilbert and Bernays, as Hilbert himself states in the preface. Other disciples, including von Neumann, were helpful but less central to the project.

methods that cannot be formalized in the system discussed.[6] But most people believed that since finitary methods are supposed to be rather weak, however they are to be delineated, they must be formalizable within any reasonably strong system. In particular, the program is supposed to be formalizable in the very system discussed. So most people concluded that Gödel had shown that Hilbert's program was hopeless. Throughout the present paper, we will ourselves assume that the consequences of the program are formalizable in the very system being discussed.

There are two relevant questions here. First, why were the Hilbertians so convinced that the program *would* work? And, second, why wasn't it eventually noted that if the program really succeeded in the way it was actually proposed in detail, it would imply its own collapse (assuming its detailed claims could be carried out in the very system discussed)?[7]

For our purposes, we will apply the Hilbert approach to systems formulated only with quantification over the natural numbers in the usual first-order logic with identity. We assume that the system —call it *S*— contains symbols for the numerals (standing for 0 and a successor function). It is also convenient to suppose that it contains function symbols for arbitrary primitive recursive functions, where the axioms imply the ordinary recursion equations. Arbitrary primitive recursive predicates can then be defined as simply predicates of the form $f(x) = 1$. Note that the proof predicate of the system, supposed to be formalized in ordinary first-order logic with identity, will itself be primitive recursive, as long as the non-logical axioms are primitive recursive. The system *S* could simply amount to the usual first-order Peano arithmetic, or it could be something stronger. The arithmetical statements provable in set theory would be much stronger, but can be separately axiomatized in accordance with these requirements using Craig's device.[8] The Hilbert approach (or *Ansatz*) was meant to apply to these stronger systems also.

---

[6] See Gödel (1931). The second incompleteness theorem is Gödel's theorem XI. In the reprinting in the Gödel papers, p. 194 (German) and p. 195 (English) contain Gödel's statement that Hilbert's formalist program is not refuted. Hilbert also stated that Gödel's work did not refute his program of the 1920s (see his introduction to Hilbert and Bernays (1934). However, as I say in the text, most people did conclude that Hilbert's original program was hopeless.

[7] I myself arrived at the present result through a circuitous route. I was looking at a purely model-theoretic formulation of the Gödel theorem (see my "A Model-Theoretic Approach to Gödel's Theorem" forthcoming) and realized that it could also be carried out syntactically, using appropriate finite approximations and semantic tableaux. But then I saw that the ladder could be kicked away and that, formulated in detail, the result the Hilbertians were attempting to obtain in fact implies its own impossibility.

[8] See Craig (1953). Of course, Craig's device would not have been known in the 1920s, but the Hilbertians could have added any extra axioms or schemata they found relevant.

The basic ideas of the program are two. Hilbert's first main idea was that quantifiers are to be eliminated, so the system will consist only of quantifier-free statements. Instead of writing $\exists x A(x)$, write $A(\varepsilon x A(x))$, where $(\varepsilon x A(x))$ denotes some true instance of $A(x)$, and is arbitrary if there is no such instance. In terms of the ε symbol, $(x)A(x)$ can be defined as $A(\varepsilon x \sim A(x))$. When all quantifiers have been eliminated in this way, one needs only the axiom scheme $A(t) \supset A(\varepsilon x A(x))$, were $t$ is any term. Terms can be formed using the constant and function symbols of the original language, but we must also allow new terms formed using the ε symbol itself. One then needs only propositional logic to deduce theorems from the axioms. Moreover, when particular values are assigned to each ε term, it is always decidable whether a given formula is true.

Then we take $A(t) \supset A(\varepsilon x A(x))$ as an axiom for the ε symbol, where $t$ is any term of the language denoting a number. $t$ can be as simple as a numeral or a primitive recursion function symbol applied to numerals, but it must be allowed to contain ε symbols itself for the idea to work out. Conditionals of this form, taken as axioms, replace the whole of quantification theory. Thus, we need only propositional logic to make inferences from the axioms. But all the axioms have to be rewritten in terms of these ε symbols. Now, stated this way, the whole thing looks pretty easy. But in fact anything written out in terms of ε symbols, as people who have studied this know, is a mess because there are quantifiers embedded inside of other quantifiers, leading to ε terms embedded in ε terms.

So, while these formulae *are* a mess, everything else is very simple because all inferences are a matter of pure propositional logic without need of quantification theory. To the Hilbertians that was important because in the intellectual atmosphere of the time (which is not today's) quantification even over all the natural numbers was—at least in classical logic, and maybe even in intuitionistic logic—a dubious idea. And Hilbert's school also officially thought that the subject of model theory (even though they themselves, in the book by Ackermann and Hilbert, raised the question of the completeness of first-order quantification theory) was meaningless.

Hilbert's second main idea was that you couldn't get a contradiction from these axioms as long as all the ε terms have true numerical values. Then, ordinary numerals for natural numbers when substituted for them will make them true. So their truth under such a substitution would show the consistency of mathematics. No contradiction could be derived from true axioms with only propositional (truth-functional) inferences.

Now, since there are no quantifiers in these formulae, whether a given substitution of natural numbers for a set of ε terms makes a formula true or false is actually decidable. Of course, they didn't have a formal theory of decidability at the time, but they knew intuitively that it was checkable. The simplest way would be to set every ε term occurring in a given formula—or throughout the proof—to zero. Then you could say whether any of the given formulae was true or false simply by deciding the zero case. Now, supposing we were lucky enough that *that* worked, we could not possibly (mistakenly) deduce a contradiction like 0 = 1, because everything is checkable in propositional logic. Of course, it is unlikely that we will be that lucky: that all the ε terms will be satisfied by the number 0. But this circumstance gave rise to the idea of another kind of interpretation, different from conventional model-theoretic interpretation; that is, a proof that is *interpreted* by giving values to the ε terms.

How can we see whether a given assignment works? Well, first we try all zeros. Probably that fails, and the formulae in the proof do not all come out true. Next, we follow a systematic procedure that the Hilbert school prescribed to allow one to change one's mind; so a 0 can change to 1 and so on.[9] Now since, as I said, ε terms may involve other ε terms, changes in the values of some ε terms intuitively ought to effect changes in the other, involved ones. On the other hand, if a given ε term was free of any dependence on other ε terms, you could change it without affecting anything else. So they set up the idea of a *priority* ordering of substitutions; that is, some changes will affect others, while other changes allow you to stick to what you already have. Now, this is more or less the idea that was rediscovered by Friedberg and Muchnik in recursion theory (computability theory): a priority ordering in which you keep on changing your mind, but which is supposed to terminate.

Hilbert and his followers tried to show that after a finite number of changes of mind, you would eventually get it right. That is, that will find values for the ε terms that will make the axioms, and hence, by propositional logic, all formulae in the proof, true. Most members of the school were convinced throughout the 1920s that this procedure just *had* to work; it is just a matter of ingenuity and detail to do it. But why were they so convinced? Well, because *really*, in the back of their minds, they thought that the axioms of first-order Peano arithmetic—or whatever the fundamental mathematical practice with numbers was —were true. And that would imply that there *are* true values for these ε terms.[10]

---

[9] I have seen at least one paper by Hilbert that formulates this in terms of the positive integers, not the non-negative integers, although this was not usual. In that case one would have to replace 0 by 1. But often he sticks with 0.

That is, if the existential statement $\exists x A(x)$ is true, then $A(\varepsilon x A(x))$, if you replace $\varepsilon x A(x)$ by some numeral, will come out true. So, therefore, they really thought it was only a matter of combinatorial work to change this into an argument that doesn't officially appeal to the model but shows that everything will terminate in values for the terms that will make the axioms true. And papers were published that purported to do this; only they contained errors. (Ackermann 1924, as described by Zach 2003, is a system of a second-order version of primitive recursion that was believed to have been proved consistent. But the consistency proof would imply the consistency of first-order Peano arithmetic, and even, if I understand matters correctly, its 1-consistency.)

The method was supposed to be completely 'finitary', but exactly what the Hilbert school meant by this is not so clear. Two clarifications have been proposed. One is Tait's— identifying finitary arguments with those that can be carried out in Skolem's primitive recursive arithmetic (Tait 1968, 1981; Skolem 1923)—but historically it appears that the Hilbert school wanted to go beyond this version, though nowadays it has been widely used as a finitary or combinatory standard. Another one is by Gödel (1958) and Kreisel (1960), who seem to believe that induction up to any fixed ordinal $< \varepsilon_0$ is finitary, but not $\varepsilon_0$ induction itself. [11]

The main requirement here is simply that the means in which the combinatorial argument is carried out be formalizable within the system itself. This requirement is much weaker than what the Hilbert school intended. It is only in virtue of our contemporary sophistication that we know that we can't reason in this way. But that way of looking at things actually influenced me, because, as I said, I started all this with model theory and only afterward translated it into proof theory. It seemed very hard at first, but finally a very simple argument was distilled.

Well, again, what shows that this program must come to an end? I was trying to say why everyone thought it *must* work. But as is well known Gödel showed that, at least if we retain the requirement that it be carried out within the system itself, it does *not*, in fact, work. Gödel's famous argument is in one sense a '*deus ex machina*.' It has no direct relation to the Hilbert program or whether it can succeed. Originally Gödel wasn't looking for such a result. His argument in the first incompleteness theorem has no direct relation to this, and he himself didn't originally realize that it was connected to the Hilbert program.[12]

---

[10] Of course, they thought the procedure would work for stronger systems also. They had a specific proposal for doing it for second order arithmetic (with quantification over number theoretic functions). But we can ignore this here. (I have little idea how, or whether, they proposed to do it for stronger systems such as *Principia*, let alone ZFC.)

[11] See Zach (2003), who states that the early paper by Ackermann (1924) used induction up to $\omega^{\omega^\omega}$, which goes beyond PRA. Zach also points out that, in this way, it resembles Gentzen's later proof of the consistency of first-order Peano arithmetic, which uses transfinite induction up to $\varepsilon_0$.

At this point let's note that the Hilbert program had, as a corollary, a much stronger result than the mere consistency of mathematics. The latter would merely mean that you couldn't prove false equations like $0 = 1$ or other false equations within the system where the identities are of arbitrary terms (not containing the ε symbol). To put it in more positive terms, any provable equations within the system will have been proved by finitary metamathematical means to be true. It will follow that any universal statement whose instances are calculable within the system must actually be true. And in most reports of the Hilbert program, either consistency is all that is stated, or more strongly that any $\Pi_1^0$ statement that is proved within the system must be true.

But the Hilbert program or approach (*Ansatz*) actually implied something much stronger. Here's the idea. Whenever a statement like $A(\varepsilon x A(x))$ is provable (remember, there are no quantifiers) there actually will be a value of the term $\varepsilon x A(x)$ that is true and provable in the system, because quantifier free calculations can be carried out within the system. This in turn implies, if we rewrite it in terms of ordinary quantification theory, that if you prove $\exists x A(x)$ and $A(x)$ is itself free of quantifiers (say, a primitive recursive predicate) some numerical instance must be provable. Therefore, we obtain something much stronger, using the details of the Hilbert program, than the mere consistency of mathematics. In contemporary consistency jargon, it would be called the 1-*consistency* of mathematics, which is a special case of what Gödel called ω-consistency (the only ω-consistency notion he actually used). But we could also call this the $\Sigma_1^0$-correctness of the system—that is, if it can prove a $\Sigma_1^0$ statement, that statement is true. Moreover, the $\Sigma_1^0$-correctness would hold in a constructive sense. One could actually calculate a true instance. Further, that von Neumann also appears to have influenced the statement of the first incompleteness theorem, at least according to Wang (1987: 43). Originally, Gödel formulated his unprovable statement as one of finite combinatorics, and von Neumann asked whether it could be made purely number-theoretic. This led Gödel to formulate his notion of Gödel numbering, and to obtain a number-theoretic independence proof.

---

[12] At least according to Dawson's biography (Dawson 1997), Gödel didn't see this right away. He received a letter from von Neumann saying: "Do you realize that it follows from your work that we can't prove mathematics consistent?" And he said, "Oh yes, I thought of that after, I lectured…" von Neumann was a nice man, at least in this respect, since others might have started raising priority issues, but he did not do so.

instance would have to be provable in the system.[13] Thus, the Hilbertian idea would have shown not merely that every $\Pi_1^0$ statement that is provable has a finitary proof, but even that the same result holds for $\Pi_2^0$ statements. For if a $\Pi_2^0$ statement is provable, each numerical instance is a provable $\Sigma_1^0$ statement, which by hypothesis must be true.[14]

The goal of the Hilbert program was to show that at least each of the narrower class of meaningful statements would have a so-called finitary proof (whatever precisely that meant)—even if that finitary proof was proved using resources that appeared to be highly *in*finitistic, i.e. resources that involved quantifiers over infinite totalities, which were therefore supposed to be suspect.

As I have said, what I am trying to show here is that although Gödel's work (independently) led to the collapse of the Hilbert program, if the Hilbertians had *really* thought through what they were claiming, they should have seen that their own demands were impossible. It is strange to me that nobody appears to have noticed this, either immediately or in the following decades. Perhaps people just regarded Gödel's work as a sufficient refutation and stopped thinking further about it. (I am struck that I myself took a surprisingly long time to see this, and only came to see it eventually via a complicated route. See footnote 7 above.)

Let's return to the question: What, therefore, is being claimed by the Hilbertians? Well, they claim that, given any proof $p$ of a $\Pi_2^0$ statement—$(x)(\exists y)A(x,y)$—where $A(x, y)$ is a simple, decidable formula not involving quantifiers, say a primitive recursive formula, and for any number $m$, there is a proof $p_1$ —remember, this can be a different proof and even longer than the proof $p$, so that $p_1$ may exceed $p$— of some particular instance $A(0^{(m)}, 0^{(n)})$. Now, first, intuitively—and reasoning infinitistically, so to speak—this is a true statement. In fact, we can eventually say something stronger than that, but at any rate this formulation is a true statement. The question is whether it is formalizable within the system. If you state this claim within the system, it says that for every proof $p$ and any number $m$ there is another proof $p_1$ and a number $n$ making an instance of this true. Notice that this statement is itself a $\Pi_2^0$ statement making a general statement about all $\Pi_2^0$ statements provable in the system.

---

[13] Of course, if a $\Sigma_1^0$ statement is true, one can calculate an instance that must be true, and therefore provable, simply by running through all the instances, together with instances of the existential quantifier until one finds one. This is a method that is at least classically valid. However, here the calculation of an instance is done more directly, and more constructively.
  Note also that the Hilbertian idea implies something even stronger: that every line of a proof containing ε terms can be replaced by numerals making these lines true. We do not need this stronger assertion here.
[14] Instead of allowing arbitrary primitive recursive predicates, one could restrict oneself to predicates involving addition, 0, S, multiplication, identity, truth-functions, and bounded quantifiers (even this is broader than is needed). But then we should stipulate that such bounded quantifiers, being decidable, do not need to be quantifiers eliminated by ε–terms. For this reason, I have chosen arbitrary primitive recursive predicates after the initial quantifiers.

The Gödel theorem has sometimes been compared to the Liar paradox. In our case we again have something like "All Cretans are liars"—only here it may seem more positive than negative. We are saying that every $\Pi_2^0$ statement has a certain property, and the statement is *itself* a $\Pi_2^0$ statement. So, whatever is self-referential here has not been introduced by some ingenious external argument—a *deus ex machina*, as I would characterize Gödel's famous argument—but rather arises from within the Hilbert program itself. They are making a $\Pi_2^0$ statement about all $\Pi_2^0$ statements, including, of course, that very statement itself.

However, we have to allow for the possibility that the Hilbertians might have realized the reflexiveness of this result while thinking of it as merely analogous to a Cretan saying, "All Cretans tell the truth". For the property they were talking about was *good*; it wasn't quite like "All Cretans are liars".

Consider the two assertions:

(*) $(x_1)(x_2)((x_1$ proves a $\Pi_2^0$ statement $\ulcorner (x)(\exists y)A(x,y) \urcorner) \supset (\exists y)(y$ proves an instance $\ulcorner A(0^{(x_2)}, 0^{(n)}) \urcorner))$

(**) $(x)((x$ proves a $\Pi_2^0$ statement $\ulcorner (x)(\exists y)A(x,y) \urcorner) \supset (\exists y)(y$ proves an instance $\ulcorner A(0^{(x)}, 0^{(n)}) \urcorner))$

Here (*) is the general claim to be made. (**) is the special case where $x_1 = x_2$.

Now (**) is a $\Pi_2^0$ statement, but to put in an appropriate normal form with one universal quantifier and one existential quantifier, we must write it as:

(***) $(x)(\exists y)((x$ proves a $\Pi_2^0$ statement $\ulcorner (x)(\exists y)A(x,y)) \urcorner \supset (y$ proves an instance $\ulcorner A(0^{(x)}, 0^{(n)}) \urcorner)$ (\*\*\*) is clearly equivalent to (\*\*). One must note here that the predicates involved are all primitive recursive. In spite of their verbal form, they need not involve quantifiers to be eliminated in the ε-calculus (since bounded quantification is primitive recursive).

Recall that we have formulated our system $S$ to be strong enough that primitive recursive predicates are correctly decidable within $S$ and are quantifier free.

(\*\*\*) is a sweeping statement about the provability of $\Pi_2^0$ statements in the system. Yet it itself is a $\Pi_2^0$ statement. Note also that it does involve the notion of Gödel numbering, assumed to be done in some standard way. Note also that the quantifier free part of (\*\*\*) is the conditional given. Call it $A^{***}(x, y)$. Now, suppose (\*\*\*) were provable. Then it has a proof with Gödel number $p$. Since $p$ does prove (\*\*\*) and true primitive recursive statements are all provable in $S$, the antecedent of $A^{***}(x, y)$ is provable, with the variable $x$ replaced by $0^{(p)}$.

Assume also that the system $S$ is 1-consistent ($\Sigma_1^0$-correct). Given this assumption, there must be a number $p_1$ instantiating the variable $y$ of $A^{***}(x, y)$. Since the antecedent of the conditional is true, the consequent must be true also. So there is a number $p_1$ that is the Gödel number of a proof of an instance of $A^{***}(x, y)$. Hence, there must be a least such $p_1$.

Now $p_1$ is supposed to be a proof of an instance of $A^{***}(0^{(p)}, y)$. But if that instance is, say, $A^{***}(0^{(p)}, 0^{(n)})$ for some particular $n$, clearly, in any standard Gödel numbering, $n < p_1$. But $n$ is supposed to be an instance of the existential quantifier in (\*\*\*). This is clearly impossible if $p_1$ is chosen to be an instance as small as possible.

(\*\*\*) is therefore unprovable if the system is 1-consistent ($\Sigma_1^0$-correct).

Now if $S$ is inconsistent, (\*\*\*) is true, since every formula is provable. On the other hand, if $S$ is 1-consistent ($\Sigma_1^0$-correct), (\*\*\*) is also true, as we have just argued. In this latter case, the situation is just like Gödel's formula: (\*\*\*) is true but unprovable.

In either case (whether $S$ is inconsistent or 1-consistent), every numerical instance of (\*\*\*) is provable, being a true $\Sigma_1^0$ statement. Hence, as in the case of Gödel's formula, (\*\*\*) cannot be refuted if $S$ is ω-consistent (in fact, 2-consistency suffices).

As in the case of Gödel's second incompleteness theorem, as long as the argument just stated can be formalized within the system $S$, it shows that '$S$ is inconsistent or 1-consistent' or equivalently, the material conditional 'if $S$ is consistent, it is 1-consistent' cannot be proved in $S$ if $S$ itself is 1-consistent. (Here, however, unlike the Gödel case, the original statement (***) comes close to being a 1-consistency statement itself.)[15]

Note that the original statement (*) also has the properties ascribed to (***). It is true if $S$ is inconsistent, and true but unprovable if $S$ is 1-consistent. Also, in either case every numerical instance (substituting for *both* variables $x_1$ and $x_2$) is provable. Here again (*) is not refutable, hence undecidable if $S$ is ω-consistent (even 2-consistent).[16] In this argument, unlike Gödel's original, the undecidability result comes from a natural statement that will be believed true by anyone who believes $S$ to be 1-consistent. The 'self-referential' aspect comes simply from the fact that it is a sweeping $\Pi_2^0$ statement about the provability of $\Pi_2^0$ statements in the system, not from an ingenious diagonalization. The statement was suggested by the Hilbert *Ansatz*, but the result and its naturalness are independent of the Hilbert program itself.[17]

*Remark 1*. What would be unknown to the Hilbertians in the argument just given? Perhaps the concept of a $\Pi_2^0$ statement, which would not have been defined in those days (but the notion seems to be implicit in their claims).

More important, the argument as given uses Gödel numbering, surely unknown to the Hilbertians. They were aware of coding devices as instantiated by Ackermann's coding of countable transfinite ordinals. But they were certainly unaware of Gödel numbering.[18] One might try to eliminate it by a direct representation of elementary syntax. However, this would require making terms cease to denote numbers and an expansion of the language.

*Remark 2*. The usual constructive proofs of the consistency of first-order Peano arithmetic using induction up to $\varepsilon_0$, starting with Gentzen (1936) and continuing with Ackermann (1940), who this time successfully used the ε-calculus, are in fact proofs of 1-consistency. They remain valid even if function symbols are added. This is the basis of Kreisel's characterization (1951, 1952) (anticipated by Gödel (1938) of the provable $\Pi_2^0$ statements of first-order Peano arithmetic, and of the $\Pi_1^1$ forms ('no counterexample interpretation') of arbitrary statements. This remains true of the later (and more elegant) proof-theoretic proofs of the consistency of first-order PA, and also the proofs (using induction up to smaller ordinals) of fragments with restricted induction.[19]

---

[15] Although, as we have seen, it is true if $S$ is inconsistent.
[16] Strictly speaking, these statements about (*) require a reformulation of (*) since it has *two* universal quantifiers, and 1-consistency, ω-consistency, etc., involve only one. But we all know how to contract quantifiers in any reasonably strong system, either using a pairing function, or more simply, a bound (i.e. $(x_1)(x_2)$ is equivalent to $(x_3)(x_1 < x_3)(x_2 < x_3)$).
[17] Of course the price is that 1-consistency needs to be assumed to get a true unprovable statement. Could one 'Rosserize' the statement to get a better result, assuming only simple consistency? I have not explored such an idea, since it would detract from the naturalness of the statement proved undecidable.
[18] In a personal communication to me, Panu Raatikainen writes: "I think the Hilbertians took for granted the possibility of something like Gödel numbering, although it was left to Gödel to actually carry out in detail one such numbering. After all, they assumed that metamathematics is done in finitary mathematics, which is essentially just a weak theory of arithmetic."

*Saul A. Kripke*
*The Saul Kripke Center and the Graduate Center of the City University of New York*

---

[19] I would like to thank Yale Weiss for editorial help. Special thanks to Romina Padró, Panu Raatikainen, and Richard Zach for their help in producing this paper. This paper has been completed with support from the Saul Kripke Center at the City University of New York, Graduate Center.